\documentclass[11pt]{article}
\usepackage{latexsym}
\usepackage{amsfonts}
\usepackage{mathrsfs,amsthm,amsmath}
\usepackage{color}
\usepackage{todonotes}
\usepackage{comment}
\newtheorem{proposition}{Proposition}[section]

\newtheorem{assumption}{Assumption}[section]

\newtheorem{lemma}{Lemma}[section]
\newtheorem{remark}{Remark}[section]
\newtheorem{example}{Example}[section]

\begin{document}
\title{Asymptotic results for sums and extremes}
\author{Rita Giuliano\thanks{Address: Dipartimento di
Matematica, Universit\`a di Pisa, Largo Bruno Pontecorvo 5,
I-56127 Pisa, Italy. e-mail: \texttt{rita.giuliano@unipi.it}} \and
Claudio Macci\thanks{Address: Dipartimento di Matematica,
Universit\`a di Roma Tor Vergata, Via della Ricerca Scientifica,
I-00133 Rome, Italy. e-mail: \texttt{macci@mat.uniroma2.it}} \and
Barbara Pacchiarotti\thanks{Address: Dipartimento di Matematica,
Universit\`a di Roma Tor Vergata, Via della Ricerca Scientifica,
I-00133 Rome, Italy. e-mail: \texttt{pacchiar@mat.uniroma2.it}}}
\date{}
\maketitle
\begin{abstract}
The term \emph{moderate deviations} is often used in the literature to mean a class of large deviation 
principles that, in some sense, fills the gap between a convergence in probability of some random variables 
to a constant, and a weak convergence to a centered Gaussian distribution (when such random variables are 
properly centered and rescaled). We talk about \emph{noncentral moderate deviations} when the weak convergence
is towards a non-Gaussian distribution. In this paper we prove a noncentral moderate deviation 
result for the bivariate sequence of sums and maxima of i.i.d. random variables bounded from above. 
We also prove a result where the random variables are not bounded from above, and the maxima are suitably normalized.
Finally we prove a moderate deviation result for sums of partial minima of i.i.d. \emph{exponential} random variables.
\end{abstract}

\section{Introduction}
The theory of large deviations gives an asymptotic computation of small probabilities on 
exponential scale; see \cite{DemboZeitouni} as a reference of this topic. The basic definition of this 
theory is the concept of \emph{large deviation principle} (LDP for short). More precisely a sequence of
probability measures $\{\pi_n:n\geq 1\}$ on a topological space $\mathcal{X}$ satisfies the LDP, with
speed $v_n$ and a rate function $I$, if the following conditions hold: $v_n\to\infty$ as $n\to\infty$,
the function $I:\mathcal{X}\to [0,\infty]$ is lower semi-continuous,
$$\liminf_{n\to\infty}\frac{1}{v_n}\log\pi_n(O)\geq-\inf_{x\in O}I(x)\ \mbox{for all open sets}\ O,$$
and
$$\limsup_{n\to\infty}\frac{1}{v_n}\log\pi_n(C)\leq-\inf_{x\in C}I(x)\ \mbox{for all closed sets}\ C.$$
Moreover we talk about weak LDP (WLDP for short) if the above upper bound for closed sets holds for compact
sets only. We also recall that, if every level set $\{x\in\mathcal{X}:I(x)\leq\eta\}$ is compact (for 
$\eta\geq 0$), the rate function $I$ is said to be \emph{good}. Finally we say that the sequence 
$\{\pi_n:n\geq 1\}$ can be seen often as a sequence of laws of $\mathcal{X}$-valued random variables 
$\{X_n:n\geq 1\}$ defined on the same probability space $(\Omega,\mathcal{F},P)$, i.e. $\pi_n=P(X_n\in\cdot)$
for every $n\geq 1$, and, with a slight abuse of terminology, we say that $\{X_n:n\geq 1\}$ satisfies the LDP.

In several common cases the rate function is good and uniquely vanishes at a certain point 
$x_\infty\in\mathcal{X}$. Then one can show that $\{X_n:n\geq 1\}$ converges in probability to $x_\infty$; 
moreover, roughly speaking, if $E$ is a Borel subset of $\mathcal{X}$ such that $x_\infty\notin\bar{E}$, 
$P(X_n\in E)$ tends to zero as $e^{-nI(E)}$ (as $n\to\infty$), where $I(E):=\inf_{x\in E}I(x)>0$.

The term \emph{moderate deviations} is used for a class of LDPs that fills the gap between two asymptotic 
regimes:
\begin{enumerate}
\item the convergence in probability of some random variables $\{X_n:n\geq 1\}$ to some constant $x_\infty$, 
governed by a LDP with speed $v_n$, and a good rate function $I$ (this LDP will called the \emph{reference LDP});
\item the weak convergence to a centered Gaussian distribution of a centered and suitably rescaled version of the 
random variables $\{X_n:n\geq 1\}$.
\end{enumerate}
More precisely we mean a class of LDPs for which the involved random variables depend on a class of positive 
sequences $\{a_n:n\geq 1\}$ (called \emph{scalings}) such that 
\begin{equation}\label{eq:MD-conditions}
a_n\to 0\ \mbox{and}\ a_nv_n\to\infty\ (\mbox{as}\ n\to\infty);
\end{equation}
the speed will be $1/a_n$ (so the speed is less fast than $v_n$ and this explains the term \emph{moderate}),
and these LDPs are governed by the same quadratic rate function $J$ which uniquely vanishes at zero. Moreover
we can say that, in some sense, one recovers the asymptotic regimes 1 and 2 above by choosing $a_n=1/v_n$ (so the 
second condition in \eqref{eq:MD-conditions} fails) and $a_n=1$ (so the first condition in \eqref{eq:MD-conditions} 
fails), respectively.

Some recent moderate deviation results in the literature concern cases in which the weak convergence is towards a
non-Gaussian distribution. Hence we talk about \emph{noncentral moderate deviations} 
(NCMD for short) and typically the common rate function $J$ for the class
of LDPs is not quadratic. Some examples of NCMD results can be found in \cite{GiulianoMacciMSTA}, where the 
weak convergences are towards Gumbel, exponential and Laplace distributions. In the same reference the interested reader 
can find references in the literature with other examples. The examples in the literature essentially concern univariate 
random variables; the only multivariate example we are aware of is the one presented in \cite{LeonenkoMacciPacchiarotti}, 
where the weak convergence is trivial because a sequence of identically distributed random variables is considered.

The aim of this paper is to prove two moderate deviation results, and a further result LDP.
\begin{itemize}
\item The first moderate deviation result fills the gap between a convergence to a constant for the bivariate
sequence of sums and maxima of i.i.d. random variables (under suitable hypotheses; in particular they are 
bounded from above), and the weak convergence towards a pair of independent random variables with standard Gaussian and Weibull 
marginal distributions (more precisely we always have the Weibull distribution of parameter 1, i.e. the distribution 
of a random variable $U$ such that $-U$ is exponentially distributed with mean 1); the weak convergence is a consequence of Theorem 1
in \cite{ChowTeugels} (see also \cite{AndersonTurkman} and \cite{Hsing} for generalizations; here we also cite 
\cite{ArendarczykKozubowskiPanorska}, \cite{Krizmanic} and \cite{QeadanKozubowskiPanorska} among the other references on the joint 
distribution of sums and maxima). Thus we obtain a NCMD result for a bivariate sequence. In particular in this paper 
we also prove the reference LDP with speed $v_n=n$, i.e. Proposition \ref{prop:LDP-sum-maxima}.
\item We prove Proposition \ref{prop:LDP-sum-maxima-added} which can be seen as a suitable modification of 
Proposition \ref{prop:LDP-sum-maxima} (under some other suitable hypotheses; in particular they are not bounded from above). As we 
shall see this new result is a LDP with $v_n=\log n$.
\item The second moderate deviation result fills the gap between a convergence to a constant for the sequence of sums of
partial minima of i.i.d. exponential random variables, and a weak convergence to a centered Gaussian distribution proved in 
\cite{Hoglund}. Thus we obtain a MD result. In this case the reference LDP with speed $v_n=\log n$ is a known
result (see Proposition 5.2 in \cite{GiulianoMacciESAIM}).
\end{itemize}
By taking into account that in this paper we present some results with sums and maxima, here we also recall some other
references which describe the impact of the maximal order statistics in the asymptotics of the sum in the normal deviation region:
\cite{Kratz}, \cite{KratzProkopenko} and \cite{Muller}. Actually these references concern the case with heavy-tailed i.i.d. random variables.

We conclude with the outline of the paper. We start with some preliminaries in Section \ref{sec:preliminaries}.
We prove the NCMD result for the bivariate sequence in Section \ref{sec:result1}, and the MD result for the sums
of partial minima of i.i.d. exponential random variables in Section \ref{sec:result2}. Finally, in 
Section \ref{sec:added}, we prove Proposition \ref{prop:LDP-sum-maxima-added}.

\section{Preliminaries}\label{sec:preliminaries}
We start with a standard way to obtain a LDP with speed $v_n$ on a Polish space $\mathcal{X}$.
Firstly, if we denote by $B_R(x)$ the open ball centered at $x$ with radius $R$, one can obtain a weak LDP (i.e. the lower bound 
for open sets, and the upper bound for compact sets) showing that
$$-I(x)\leq\lim_{R\to 0}\liminf_{n\to\infty}\frac{1}{v_n}\log\pi_n(B_R(x))\leq
\lim_{R\to 0}\limsup_{n\to\infty}\frac{1}{v_n}\log\pi_n(B_R(x))\leq-I(x)$$
for all $x\in\mathcal{X}$ (this can be seen as a consequence of Theorem 4.1.11 in \cite{DemboZeitouni}); actually one can consider an 
arbitrary basis of neighborhoods of each point $x\in\mathcal{X}$ instead of open balls. Successively (see e.g. Lemma 1.2.18 in 
\cite{DemboZeitouni}) one can obtain the full LDP (i.e. the upper bound for closed sets) showing that the \emph{exponential tightness} 
condition holds (see e.g. \cite{DemboZeitouni}, page 8): for all $b>0$ there exists a compact set $K_b$ such that
$$\limsup_{n\to\infty}\frac{1}{v_n}\log\pi_n(K_b^c)\leq -b$$
or, equivalently,
$$\pi_n(K_b^c)\leq ae^{-v_nb}\quad \mbox{eventually}$$
for some $a>0$.

Moreover, in view of what we present in the next sections, we recall two results. The first one is the well-known 
G\"artner Ellis Theorem (see e.g. Theorem 2.3.6(c) in \cite{DemboZeitouni}). Here we recall a simplified version of the theorem, with 
$\mathcal{X}=\mathbb{R}$.

\begin{proposition}\label{prop:GET}
Let $\{\pi_n:n\geq 1\}$ be a sequence of probability measures on $\mathbb{R}$, and let $\{v_n:n\geq 1\}$ be a speed function. 
Moreover assume that, for all $\theta\in\mathbb{R}$, the limit
$$\Lambda(\theta):=\lim_{n\to\infty}\frac{1}{v_n}\int_{\mathbb{R}}e^{v_n\theta x}\pi_n(dx)$$
exists as an extended real number, and that $0\in(\mathcal{D}(\Lambda))^\circ$, where 
$\mathcal{D}(\Lambda):=\{\theta\in\mathbb{R}:\Lambda(\theta)<\infty\}$. Then, if $\Lambda$ is essentially smooth and lower 
semi-continuous, the sequence $\{\pi_n:n\geq 1\}$ satisfies the LDP with good rate function $\Lambda^*$ defined by 
$\Lambda^*(x):=\sup_{\theta\in\mathbb{R}}\{\theta x-\Lambda(\theta)\}$.
\end{proposition}

For completeness we recall that the function $\Lambda$ is essentially smooth (see e.g. Definition 2.3.5 in \cite{DemboZeitouni})
if it is differentiable throughout the set $(\mathcal{D}(\Lambda))^\circ$ (assumed to be non-empty), and if it is steep, i.e. 
$|\Lambda^\prime(\theta_n)|$ tends to infinity whenever $\{\theta_n:n\geq 1\}\subset(\mathcal{D}(\Lambda))^\circ$ is a sequence 
which converges to a boundary point of $\mathcal{D}(\Lambda)$.

The second result is Theorem 2.3 in \cite{Chaganty} which plays a crucial role in the proofs of the first moderate 
deviation result (i.e. the one for the bivariate sequence of sums and maxima of i.i.d. random variables). In this case we have 
$\mathcal{X}=\mathcal{Y}\times\mathcal{Z}$, where $\mathcal{Y}$ and $\mathcal{Z}$ are Polish spaces. Moreover, for any given
probability measure $\pi_n$ on $\mathcal{X}=\mathcal{Y}\times\mathcal{Z}$, we consider the marginal distributions
$\pi_n^Y$ on $\mathcal{Y}$ and $\pi_n^Z$ on $\mathcal{Z}$, i.e.
$$\pi_n^Y(dy)=\int_{\mathcal{Z}}\pi_n(dy,dz)\ \mbox{and}\ \pi_n^Z(dz)=\int_{\mathcal{Y}}\pi_n(dy,dz),$$
and the conditional distributions $\{\pi_n^{Y|Z}(dy|z):z\in\mathcal{Z}\}$ on $\mathcal{Y}$ such that
$$\pi_n(dy,dz)=\pi_n^{Y|Z}(dy|z)\pi_n^Z(dz).$$

\begin{proposition}\label{prop:Chaganty}
Let $\{\pi_n:n\geq 1\}$ be a sequence of probability measures on $\mathcal{X}=\mathcal{Y}\times\mathcal{Z}$, where $\mathcal{Y}$ 
and $\mathcal{Z}$ are Polish spaces. We assume that the following conditions hold.\\
$\mathbf{(C1)}$: The sequence $\{\pi_n^Z:n\geq 1\}$ satisfies the LDP with speed $v_n$ and a good rate function $I_Z$.\\
$\mathbf{(C2)}$: If $\{z_n:n\geq 1\}\subset\mathcal{Z}$ and $z_n\to z\in\mathcal{Z}$, then $\{\pi_n^{Y|Z}(dy|z_n):n\geq 1\}$
satisfies the LDP with speed $v_n$ and good rate function $I_{Y|Z}(\cdot|z)$, where $\{I_{Y|Z}(\cdot|z):z\in\mathcal{Z}\}$
is a family of good rate functions such that
\begin{equation}\label{eq:joint-lsc}
(y,z)\mapsto I_{Y|Z}(y|z)\ \mbox{is lower semicontinuous}.
\end{equation}
Then $\{\pi_n:n\geq 1\}$ satisfies the WLDP with speed $v_n$ and rate function $I_{Y,Z}$ defined by
$$I_{Y,Z}(y,z):=I_{Y|Z}(y|z)+I_Z(z).$$
Moreover:
$\{\pi_n^Y:n\geq 1\}$ satisfies the LDP with speed $v_n$ and rate function $I_Y$ defined by
$$I_Y(y):=\inf_{z\in\mathcal{Z}}\{I_{Y,Z}(y,z)\}=\inf_{z\in\mathcal{Z}}\{I_{Y|Z}(y|z)+I_Z(z)\};$$
$\{\pi_n:n\geq 1\}$ satisfies the full LDP if the rate function $I_{Y,Z}$ is good and, in such a case, the rate function $I_Y$ 
is also good.
\end{proposition}

In what follows we apply Proposition \ref{prop:Chaganty} in the proofs of Propositions \ref{prop:LDP-sum-maxima} and 
\ref{prop:sum-maxima-NCMD}. Actually we omit the statement for $\{\pi_n^Y:n\geq 1\}$ because it would allow to recover 
well-known results.

\section{NCMD for sums and maxima of i.i.d. random variables}\label{sec:result1}
Throughout this section we assume the following

\begin{assumption}\label{assumpt:LD}
Let $\{W_n:n\geq 1\}$ be a sequence of i.i.d. real random variables with density function $f$ which is
assumed to be positive only on an interval $(m,M)$, where $-\infty\leq m<M<+\infty$. We set
$$\mathcal{I}=\overline{(m,M)}=\left\{\begin{array}{ll}
[m,M]&\ \mbox{if}\ m>-\infty\\
(-\infty,M]&\ \mbox{if}\ m=-\infty.
\end{array}\right.$$
Moreover, as usual, we set $F(z):=\int_{-\infty}^zf(w)dw$ for $z\in\mathbb{R}$; then $F(M)=1$ and, if 
$m>-\infty$, $F(m)=0$. Finally we also assume 
that, for every $z\in\mathcal{I}$, the function $\kappa_{Y|Z}(\cdot|z)$ defined by
\begin{equation}\label{eq:def-kappa-Y|Z}
\kappa_{Y|Z}(\theta|z):=\left\{\begin{array}{ll}
\theta m&\ \mbox{if}\ z=m>-\infty\\
\log\frac{\int_{-\infty}^z e^{\theta w}f(w)dw}{F(z)}&\ \mbox{otherwise}.
\end{array}\right.
\end{equation}
is finite in a neighborhood of the origin $\theta=0$, essentially smooth and lower semicontinuous.
\end{assumption}

We are interested in the asymptotic behavior of the sequence of bivariate random variables 
$\{(Y_n,Z_n):n\geq 1\}$ defined by
$$(Y_n,Z_n):=\left(\frac{W_1+\cdots+W_n}{n},\max\{W_1,\ldots,W_n\}\right).$$
The first result in this section is the reference LDP, i.e. Proposition \ref{prop:LDP-sum-maxima} which provides
the full LDP of $\{P((Y_n,Z_n)\in\cdot):n\geq 1\}$ in the final part of the statement of the proposition. 
Successively, under some further conditions (see Assumption \ref{assumpt:NCMD} below), in Proposition 
\ref{prop:ChowTeugels-consequence} we show that we have a weak convergence towards a non-Gaussian distribution, 
and in Proposition \ref{prop:sum-maxima-NCMD} we prove a NCMD result. Both Propositions 
\ref{prop:LDP-sum-maxima} and \ref{prop:sum-maxima-NCMD} will be proved by applying Theorem 2.3 in \cite{Chaganty},
i.e. Proposition \ref{prop:Chaganty}.

\subsection{The reference LDP}
We start with the following

\begin{proposition}\label{prop:LDP-sum-maxima}
Assume that Assumption \ref{assumpt:LD} holds. Let $I_Z$ be defined by 
$$I_Z(z):=-\log F(z),\ \mbox{for}\ z\in\mathcal{I},$$
(with the rule $\log 0=-\infty$ for the case $z=m$ when $m>-\infty$), and let $\{I_{Y|Z}(\cdot|z):z\in\mathcal{I}\}$ 
be the functions defined by
$$I_{Y|Z}(y|z):=\sup_{\theta\in\mathbb{R}}\{\theta y-\kappa_{Y|Z}(\theta|z)\},$$
where $\kappa_{Y|Z}(\theta|z)$ is defined by \eqref{eq:def-kappa-Y|Z}. Then $\{P((Y_n,Z_n)\in\cdot):n\geq 1\}$
satisfies the WLDP with speed $n$ and rate function $I_{Y,Z}$ defined by
$$I_{Y,Z}(y,z):=I_{Y|Z}(y|z)+I_Z(z).$$
Moreover $\{P((Y_n,Z_n)\in\cdot):n\geq 1\}$ satisfies the full LDP if the rate function $I_{Y,Z}$ is good 
and, in such a case, the rate function $I_Y$ is also good.
\end{proposition}
\begin{proof}
We want to apply Proposition \ref{prop:Chaganty} (on the product space 
$\mathcal{Y}\times\mathcal{Z}:=\mathcal{I}\times\mathcal{I}$) to the sequence $\{\pi_n:n\geq 1\}$ defined by
$$\pi_n(\cdot)=P((Y_n,Z_n)\in\cdot).$$
It is known that Condition $\mathbf{(C1)}$ trivially holds; see e.g. Proposition 4.1 in \cite{GiulianoMacciCSTM}. 
So, in the remaining part of the proof, we have to show that Condition $\mathbf{(C2)}$ holds.

Firstly we can easily check condition \eqref{eq:joint-lsc}. Indeed, if we take
$\{(y_n,z_n):n\geq 1\}\subset\mathcal{I}\times\mathcal{I}$ such that $(y_n,z_n)\to (y,z)\in\mathcal{I}\times\mathcal{I}$,
we have
$$I_{Y|Z}(y_n|z_n)\geq\theta y_n-\kappa_{Y|Z}(\theta|z_n)\ \mbox{for all}\ \theta\in\mathbb{R},$$
which yields (if $z>m$ this is trivial; if $z=m>-\infty$ this follows from an application of L'H\^{o}pital's rule)
$$\liminf_{n\to\infty}I_{Y|Z}(y_n|z_n)\geq\theta y-\kappa_{Y|Z}(\theta|z)\ \mbox{for all}\ \theta\in\mathbb{R},$$
and we get
$$\liminf_{n\to\infty}I_{Y|Z}(y_n|z_n)\geq I_{Y|Z}(y|z)$$
by taking the supremum with respect to $\theta\in\mathbb{R}$. Thus \eqref{eq:joint-lsc} is checked.

In order to complete the proof of Condition $\mathbf{(C2)}$, some preliminaries are needed. Namely we recall a well-known 
result on order statistics, and we introduce a suitable family of densities.
\begin{itemize}
\item For every $n\geq 1$, let $W_{1:n},\ldots,W_{n:n}$ be the order statistics of $W_1,\ldots,W_n$.
Then the joint distribution of $(W_{1:n},\ldots,W_{n:n})$ has density
$$g(w_1,\ldots,w_n)=n!f(w_1)\cdots f(w_n)1_{w_1<\cdots<w_n}.$$
\item For every $z\in\mathcal{I}$ such that $z\neq m$ when $m>-\infty$, we introduce the density $f(\cdot|z)$ defined by
\begin{equation}\label{eq:restricted-density}
f(w|z)=\frac{f(w)}{F(z)}1_{(-\infty,z)}(w).
\end{equation}
\end{itemize}
We assume for the moment that
\begin{equation}\label{eq:Darling-replace}
\log\mathbb{E}[e^{n\theta Y_n}|Z_n=z_n]=(n-1)\kappa_{Y|Z}(\theta|z_n)+\theta z_n\ (\mbox{for all}\ \theta\in\mathbb{R});
\end{equation}
this will be checked below. Then, by \eqref{eq:Darling-replace}, we get
\begin{equation}\label{eq:GET-condition}
\lim_{n\to\infty}\frac{1}{n}\log\mathbb{E}[e^{n\theta Y_n}|Z_n=z_n]=\kappa_{Y|Z}(\theta|z)\ (\mbox{for all}\ \theta\in\mathbb{R});
\end{equation}
and, by the hypotheses on the functions $\{\kappa_{Y|Z}(\cdot|z):z\in\mathcal{I}\}$, we see that Condition $\mathbf{(C2)}$ holds
by a straightforward application of G\"artner Ellis Theorem, i.e. Proposition \ref{prop:GET}.

In conclusion we have to check \eqref{eq:Darling-replace}, and for simplicity we write $z$ in place of $z_n$. Actually the case 
$z=m$, when $m>-\infty$, is immediate; therefore, from now on, we assume to have $z>m$. Firstly we have
\begin{multline*}
\mathbb{E}[e^{n\theta Y_n}|Z_n=z]=\mathbb{E}[e^{\theta\sum_{i=1}^nW_{i:n}}|W_{n:n}=z]=
\int_{\mathbb{R}^{n-1}}e^{\theta(\sum_{i=1}^{n-1}w_i+z)}\ \frac{g(w_1,\ldots,w_{n-1},z)}{n(F(z))^{n-1}f(z)}dw_1\cdots dw_{n-1}\\
=e^{\theta z}\int_{\mathbb{R}^{n-1}}e^{\theta\sum_{i=1}^{n-1}w_i}\ \frac{n!f(w_1)\cdots f(w_{n-1})f(z)1_{w_1<\cdots<w_{n-1}<z}}
{n(F(z))^{n-1}f(z)}dw_1\cdots dw_{n-1}\\
=e^{\theta z}\int_{\mathbb{R}^{n-1}}e^{\theta\sum_{i=1}^{n-1}w_i}\ \frac{(n-1)!f(w_1)\cdots f(w_{n-1})1_{w_1<\cdots<w_{n-1}<z}}
{(F(z))^{n-1}}dw_1\cdots dw_{n-1}\\
=e^{\theta z}\int_{\mathbb{R}^{n-1}}e^{\theta\sum_{i=1}^{n-1}w_i}\ (n-1)!f(w_1|z)\cdots f(w_{n-1}|z)1_{w_1<\cdots<w_{n-1}}dw_1\cdots dw_{n-1}.
\end{multline*}
Then let $W_1^{(z)},\ldots,W_{n-1}^{(z)}$ be i.i.d. random variables with common density $f(\cdot|z)$ in \eqref{eq:restricted-density};
moreover we denote their order statistics by $W_{1:n-1}^{(z)},\ldots,W_{n-1:n-1}^{(z)}$. We have
$$\mathbb{E}[e^{n\theta Y_n}|Z_n=z]=e^{\theta z}\mathbb{E}[e^{\theta\sum_{i=1}^{n-1}W_{i:n-1}^{(z)}}]
=e^{\theta z}\mathbb{E}[e^{\theta\sum_{i=1}^{n-1}W_i^{(z)}}]
=e^{\theta z}\prod_{i=1}^{n-1}\mathbb{E}[e^{\theta W_i^{(z)}}]$$
and, by taking into account the definition of the function $\kappa_{Y|Z}(\cdot|z)$ in \eqref{eq:def-kappa-Y|Z}, we get
$$\mathbb{E}[e^{n\theta Y_n}|Z_n=z]=e^{\theta z}(e^{\kappa_{Y|Z}(\theta|z)})^{n-1}
=e^{(n-1)\kappa_{Y|Z}(\theta|z)+\theta z}.$$
Thus \eqref{eq:Darling-replace} is checked.
\end{proof}

We have the following remarks.

\begin{remark}\label{rem:compact-case}
The sequence $\{P((Y_n,Z_n)\in\cdot):n\geq 1\}$ satisfies the full LDP (because the rate function $I_{Y,Z}$ is good) if $m>-\infty$, 
i.e. if $\mathcal{I}$ is compact. Indeed, in such a case, every closed level set of $I_{Y,Z}$ is compact (indeed it is a subset of 
$\mathcal{I}\times\mathcal{I}$). We also recall that, if $m>-\infty$, the function $\kappa_{Y|Z}(\cdot|z)$ is finite and differentiable.
\end{remark}

\begin{remark}\label{rem:beyond-M<infty}
One can wonder if we can obtain a version of Proposition \ref{prop:LDP-sum-maxima} when the distribution of the random variables 
$\{W_n:n\geq 1\}$ is not bounded from above, i.e. when $M=\infty$. Firstly, in such a case, the rate function $I_Z$ in Proposition 
\ref{prop:LDP-sum-maxima} is not good, and therefore we cannot apply Proposition \ref{prop:Chaganty}. However we can prove Proposition 
\ref{prop:LDP-sum-maxima-added}, i.e. a suitable modification of Proposition \ref{prop:LDP-sum-maxima} with 
$\{P((Y_n,Z_n/h_n)\in\cdot):n\geq 1\}$ in place of $\{P((Y_n,Z_n)\in\cdot):n\geq 1\}$, for some $h_n$ such that $h_n\to\infty$.
\end{remark}

It is interesting to present the following example in which the rate function $I_{Y,Z}$ is good even if $m=-\infty$.

\begin{example}\label{ex:non-compact}
We take $\mathcal{I}=(-\infty,0]$ (so $M=0$ and $m=-\infty$). Let $f$ be 
defined by
$$f(w):=e^w1_{(-\infty,0)}(w).$$
Then, for $z\in (-\infty,0]$, we have $F(z)=e^z$ which yields: $I_Z(z)=-z$;
\begin{multline*}
\kappa_{Y|Z}(\theta|z)=\log\frac{\int_{-\infty}^ze^{\theta w}e^wdw}{e^z}=\log\left(e^{-z}\int_{-\infty}^ze^{(\theta+1)w}dw\right)\\
=\left\{\begin{array}{ll}
-z+\log\frac{e^{(\theta+1)z}}{\theta+1}&\ \mbox{if}\ \theta+1>0\\
\infty&\ \mbox{if}\ \theta+1\leq 0
\end{array}\right.=\left\{\begin{array}{ll}
\theta z-\log(\theta+1)&\ \mbox{if}\ \theta>-1\\
\infty&\ \mbox{if}\ \theta\leq -1,
\end{array}\right.
\end{multline*}
and therefore, for $y<z\leq 0$,
$$I_{Y|Z}(y|z)=\sup_{\theta\in\mathbb{R}}\{\theta y-\kappa_{Y|Z}(\theta|z)\}
=\sup_{\theta>-1}\{\theta(y-z)+\log(\theta+1)\}=(z-y)-1-\log(z-y)$$
(indeed the supremum above is attained at $\theta=\frac{1}{z-y}-1$).

Now we are able to check that $I_{Y,Z}$ is a good rate function. For every $\eta\geq 0$ we have
\begin{multline*}
\{(y,z)\in\mathcal{I}\times\mathcal{I}:I_{Y,Z}(y,z)\leq\eta\}=\{(y,z)\in\mathcal{I}\times\mathcal{I}:I_{Y|Z}(y|z)+I_Z(z)\leq\eta\}\\
=\{(y,z):y<z\leq 0,(z-y)-1-\log(z-y)-z\leq\eta\}\\
\subset\{(y,z):y<z\leq 0,(z-y)-1-\log(z-y)\leq\eta,-z\leq\eta\};
\end{multline*}
moreover, for every $z\in[-\eta,0]$, there exist $r_1^\eta,r_2^\eta$ such that
$0<r_1^\eta<1<r_2^\eta$ and
$$(z-y)-1-\log(z-y)\leq\eta,\ \mbox{if and only if}\ r_1^\eta<z-y<r_2^\eta;$$
therefore the (closed) level set $\{(y,z)\in\mathcal{I}\times\mathcal{I}:I_{Y,Z}(y,z)\leq\eta\}$ is a subset of
the compact set (it is a parallelogram)
$$\{(y,z):-\eta\leq z\leq 0,z-r_2^\eta\leq y\leq z-r_1^\eta\}.$$
\end{example}

\subsection{The weak convergence and NCMD}
Throughout this paper we consider the Weibull distribution with parameter 1, i.e. the distribution of a random variable $U$ such that
$$P(U\leq u)=\min\{e^u,1\}\ \mbox{for all}\ u\in\mathbb{R}$$
(thus $-U$ is an exponentially distributed random variable with mean 1). We start with the following assumption.

\begin{assumption}\label{assumpt:NCMD}
Let $\{W_n:n\geq 1\}$ be a sequence of i.i.d. real random variables as in Assumption \ref{assumpt:LD} with density 
function $f$ (so, in particular, the random variables $\{W_n:n\geq 1\}$ have finite mean $\mu<M$ and variance 
$\sigma^2>0$, indeed $\kappa_{Y|Z}(\theta|M)=\log\mathbb{E}[e^{\theta W_1}]$ is finite in a neighborhood of the origin $\theta=0$).
Moreover we assume that there exists $F^\prime(M-)>0$, i.e. the left derivative of $F$ at $M$, and that $f(M)=F^\prime(M-)$. Finally let 
$\{L(n):n\geq 1\}$ be a sequence such that $L(n)\to F^\prime(M-)$ as $n\to\infty$. 
\end{assumption}

It is well-known that, if Assumption \ref{assumpt:NCMD} holds, we have the following weak convergence results 
(as $n\to\infty$):
\begin{enumerate}
\item $\left\{\frac{Y_n-\mu}{\sigma/\sqrt{n}}:n\geq 1\right\}$ converges weakly to a standard Gaussian distribution 
(by the Central Limit Theorem);
\item $\{nL(n)(Z_n-M):n\geq 1\}$ converges weakly to a Weibull distribution with parameter 1. Indeed, for every $z\leq 0$,
for a suitable remainder $o\left(\frac{1}{n}\right)$ (as $n\to\infty$), for $n$ large enough we have
\begin{multline*}
P(nL(n)(Z_n-M)\leq z)=P\left(Z_n\leq M+\frac{z}{nL(n)}\right)\\
=F^n\left(M+\frac{z}{nL(n)}\right)=\left(1+F^\prime(M-)\frac{z}{nL(n)}+o\left(\frac{1}{n}\right)\right)^n
\to e^z\ (\mbox{as}\ n\to\infty).
\end{multline*}
\end{enumerate}

\begin{remark}\label{rem:EVT-connection}
The weak convergence of $\{nL(n)(Z_n-M):n\geq 1\}$ in item 2 above can be related to a particular case of the Fisher 
Tippett Gnedenko Theorem (see e.g. Theorem 3.2.3 in \cite{EmbrechtsKluppelbergMikosch}). More precisely we mean the weak convergence
of $\left\{\frac{Z_n-M}{M-F^{-1}(1-1/n)}:n\geq 1\right\}$ when the random variables $\{W_n:n\geq 1\}$ are in the 
\emph{Maximum Domain of Attraction} of the Weibull distribution with parameter 1; see e.g. the summarizing Table 3.4.3 in 
\cite{EmbrechtsKluppelbergMikosch}, page 154, for $\alpha=1$ (for the related theory see Section 3.3.2 in the same reference).
Indeed we have
$$M-F^{-1}(1-1/n)=n^{-1}L_1(n)$$
for a slowly varying function $L_1$; then, since $M=F^{-1}(1)$, we get
$$L_1(n)=\frac{F^{-1}(1)-F^{-1}(1-1/n)}{1/n}\to (F^{-1})^\prime(1)=\frac{1}{F^\prime(M-)}\ \mbox{as}\ n\to\infty,$$
and therefore $L_1(n)$ plays the role of $\frac{1}{L(n)}$ (at least for $n$ large enough).
\end{remark}

Actually, as we say in the next proposition, we have a stronger result, i.e. the weak convergence of the bivariate 
sequence towards a bivariate distribution with independent components.

\begin{proposition}\label{prop:ChowTeugels-consequence}
If Assumption \ref{assumpt:NCMD} holds, then $\left\{\left(\frac{Y_n-\mu}{\sigma/\sqrt{n}},nL(n)(Z_n-M)\right):n\geq 1\right\}$
converges weakly to a bivariate distribution with independent components distributed as a standard Gaussian distribution and a 
Weibull distribution with parameter 1.
\end{proposition}
\begin{proof}
It is a consequence of Theorem 1 in \cite{ChowTeugels}.
\end{proof}

The next Proposition \ref{prop:sum-maxima-NCMD} provides a class of LDPs that fills the gap between the convergence of 
$\{(Y_n,Z_n):n\geq 1\}$ to $(\mu,M)$ (governed by the LDP in Proposition \ref{prop:LDP-sum-maxima} with speed $v_n=n$), and 
the weak convergence in Proposition \ref{prop:ChowTeugels-consequence}. Then we have a NCMD result 
because the weak convergence in Proposition \ref{prop:ChowTeugels-consequence} is towards a non-Gaussian distribution
(indeed the second marginal distribution is not Gaussian).

\begin{proposition}\label{prop:sum-maxima-NCMD}
Assume that Assumption \ref{assumpt:NCMD} holds. Then, for every sequence of positive numbers $\{a_n:n\geq 1\}$ such that
\eqref{eq:MD-conditions} holds with $v_n=n$, 
$\left\{P\left(\left(\sqrt{a_n}\frac{Y_n-\mu}{\sigma/\sqrt{n}},a_nnL(n)(Z_n-M)\right)\in\cdot\right):n\geq 1\right\}$
satisfies the LDP with speed $1/a_n$ and good rate function $J_{Y,Z}$ defined by
$$J_{Y,Z}(y,z)=\left\{\begin{array}{ll}
\frac{y^2}{2}-z&\ \mbox{if}\ z\leq 0\\
\infty&\ \mbox{otherwise}.
\end{array}\right.$$
\end{proposition}
\begin{proof}
We want to apply Proposition \ref{prop:Chaganty} (on the product space $\mathcal{Y}\times\mathcal{Z}:=\mathbb{R}\times (-\infty,0]$)
to the sequence $\{\pi_n:n\geq 1\}$ defined by
$$\pi_n(\cdot)=P\left(\left(\sqrt{a_n}\frac{Y_n-\mu}{\sigma/\sqrt{n}},a_nnL(n)(Z_n-M)\right)\in\cdot\right).$$
Notice that here we use some slightly different notation (i.e. $J_Z$, $J_{Y|Z}$ and $J_{Y,Z}$ in place of $I_Z$, $I_{Y|Z}$
and $I_{Y,Z}$ in Proposition \ref{prop:Chaganty}, respectively). The proof is divided in three parts: the proof of Condition 
$\mathbf{(C1)}$, the proof of Condition $\mathbf{(C2)}$, and the proof of the goodness of the rate function $J_{Y,Z}$.

\paragraph{The proof of Condition $\mathbf{(C1)}$.}
Here we consider the sequence $\{\pi_n^Z:n\geq 1\}$ defined by
$$\pi_n^Z(\cdot)=P(a_nnL(n)(Z_n-M)\in\cdot).$$
Then we have to prove that $\{\pi_n^Z:n\geq 1\}$ satisfies the LDP with speed $1/a_n$ and good rate function $J_Z$ 
defined by
$$J_Z(z)=\left\{\begin{array}{ll}
-z&\ \mbox{if}\ z\leq 0\\
\infty&\ \mbox{otherwise}.
\end{array}\right.$$
We start with the proof of the upper bound for every closed set $C\subset(-\infty,0]$. If $0\in C$ it is trivial. If $0\notin C$, 
we set $z_C:=\sup C$ and therefore we have
$$z_C=-\inf_{z\in C} J_Z(z)<0,\ \mbox{with}\ z_C\in C.$$
Then, for a suitable remainder $o\left(\frac{1}{a_nn}\right)$ (as $n\to\infty$), for $n$ large enough we have
\begin{multline*}
P(a_nnL(n)(Z_n-M)\in C)\leq P(a_nnL(n)(Z_n-M)\leq z_C)\\
=P\left(Z_n\leq M+\frac{z_C}{a_nnL(n)}\right)
=F^n\left(M+\frac{z_C}{a_nnL(n)}\right)=\left(1+F^\prime(M-)\frac{z_C}{a_nnL(n)}+o\left(\frac{1}{a_nn}\right)\right)^n,
\end{multline*}
and therefore
\begin{multline*}
\limsup_{n\to\infty}\frac{1}{1/a_n}\log P(a_nnL(n)(Z_n-M)\in C)\\
\leq\limsup_{n\to\infty}a_nn\log\left(1+F^\prime(M-)\frac{z_C}{a_nnL(n)}+o\left(\frac{1}{a_nn}\right)\right)
=z_C=-\inf_{z\in C} J_Z(z).
\end{multline*}
Now the lower bound for open sets. For every open set $O\in(-\infty,0]$ such that $z\in O$, we have to check that
$$\limsup_{n\to\infty}\frac{1}{1/a_n}\log P(a_nnL(n)(Z_n-M)\in O)\geq-J_Z(z).$$
This is trivial if $z=0$ because $P(a_nnL(n)(Z_n-M)\in O)\to 1$ because $a_nnL(n)(Z_n-M)$ converges 
in probability to zero (it is a trivial consequence of the Slutsky Theorem). For $z<0$ we take $\varepsilon>0$ small
enough to have $(z-\varepsilon,z+\varepsilon)\subset O\cap(-\infty,0)$ and, by also taking into account some 
computations above for the proof of the upper bound for closed sets, for $n$ large enough we get
\begin{multline*}
P(a_nnL(n)(Z_n-M)\in O)\geq P(z-\varepsilon<a_nnL(n)(Z_n-M)<z+\varepsilon)\\
=P\left(M+\frac{z-\varepsilon}{a_nnL(n)}<Z_n<M+\frac{z+\varepsilon}{a_nnL(n)}\right)\\
=F^n\left(M+\frac{z+\varepsilon}{a_nnL(n)}\right)-F^n\left(M+\frac{z-\varepsilon}{a_nnL(n)}\right)\\
=\left(1+F^\prime(M-)\frac{z+\varepsilon}{a_nnL(n)}+o\left(\frac{1}{a_nn}\right)\right)^n
-\left(1+F^\prime(M-)\frac{z-\varepsilon}{a_nnL(n)}+o\left(\frac{1}{a_nn}\right)\right)^n\\
=\left(1+F^\prime(M-)\frac{z-\varepsilon}{a_nnL(n)}+o\left(\frac{1}{a_nn}\right)\right)^n
\left(\frac{\left(1+F^\prime(M-)\frac{z+\varepsilon}{a_nnL(n)}+o\left(\frac{1}{a_nn}\right)\right)^n}
{\left(1+F^\prime(M-)\frac{z-\varepsilon}{a_nnL(n)}+o\left(\frac{1}{a_nn}\right)\right)^n}-1\right);
\end{multline*}
moreover
\begin{multline*}
\liminf_{n\to\infty}\frac{1}{1/a_n}\log	P(a_nnL(n)(Z_n-M)\in O)\geq
\liminf_{n\to\infty}a_nn\log\left(1+F^\prime(M-)\frac{z-\varepsilon}{a_nnL(n)}+o\left(\frac{1}{a_nn}\right)\right)\\
+\liminf_{n\to\infty}a_n\log\left(\exp\left(n\log\left(\frac{1+F^\prime(M-)\frac{z+\varepsilon}{a_nnL(n)}+o\left(\frac{1}{a_nn}\right)}
{1+F^\prime(M-)\frac{z-\varepsilon}{a_nnL(n)}+o\left(\frac{1}{a_nn}\right)}\right)\right)-1\right),
\end{multline*}
where
$$\liminf_{n\to\infty}a_nn\log\left(1+F^\prime(M-)\frac{z-\varepsilon}{a_nnL(n)}+o\left(\frac{1}{a_nn}\right)\right)=z-\varepsilon$$
and
$$n\log\left(\frac{1+F^\prime(M-)\frac{z+\varepsilon}{a_nnL(n)}+o\left(\frac{1}{a_nn}\right)}
{1+F^\prime(M-)\frac{z-\varepsilon}{a_nnL(n)}+o\left(\frac{1}{a_nn}\right)}\right)
=n\log\left(1+\frac{F^\prime(M-)\frac{2\varepsilon}{a_nnL(n)}+o\left(\frac{1}{a_nn}\right)}
{1+F^\prime(M-)\frac{z-\varepsilon}{a_nnL(n)}+o\left(\frac{1}{a_nn}\right)}\right)
\sim\frac{2\varepsilon}{a_n};$$
so finally we have
$$\liminf_{n\to\infty}\frac{1}{1/a_n}\log P(a_nnL(n)(Z_n-M)\in O)\geq z-\varepsilon+2\varepsilon=-J_Z(z)+\varepsilon,$$
and we conclude by letting $\varepsilon$ go to zero.

\paragraph{The proof of Condition $\mathbf{(C2)}$.}
Here we consider the sequence $\{\pi_n^{Y|Z}(\cdot|z_n):n\geq 1\}$ defined by
$$\pi_n^{Y|Z}(\cdot|z_n)=P\left(\sqrt{a_n}\frac{Y_n-\mu}{\sigma/\sqrt{n}}\in\cdot\Big|a_nnL(n)(Z_n-M)=z_n\right),$$
where $\{z_n:n\geq 1\}\subset (-\infty,0]$ such that $z_n\to z$ (as $n\to\infty$) for some $z\in (-\infty,0]$.
Then we have to prove that $\{\pi_n^{Y|Z}(\cdot|z_n):n\geq 1\}$ satisfies the LDP with speed $1/a_n$ and good rate function $J_{Y|Z}$ 
defined by
$$J_{Y|Z}(y|z)=\frac{y^2}{2}.$$
Note that condition \eqref{eq:joint-lsc} trivially holds; indeed $(y,z)\mapsto J_{Y|Z}(y|z)=\frac{y^2}{2}$ is a lower semicontinuous 
function. Moreover, in what follows, we simply write $J_Y(y)=\frac{y^2}{2}$ in place of $J_{Y|Z}(y|z)=\frac{y^2}{2}$.

We apply G\"artner Ellis Theorem, i.e. Proposition \ref{prop:GET}. Indeed we show that
\begin{equation}\label{eq:GETlimit-conditionals-NCMD}
\lim_{n\to\infty}\frac{1}{1/a_n}\log\mathbb{E}\left[\exp\left(\frac{\theta}{a_n}\sqrt{a_n}\frac{Y_n-\mu}{\sigma/\sqrt{n}}\right)
\Big|a_nnL(n)(Z_n-M)=z_n\right]=\frac{\theta^2}{2}\ (\mbox{for all}\ \theta\in\mathbb{R})
\end{equation}
and therefore, for every $z\leq 0$, we get the desired LDP with rate function $J_Y$ defined by
$$J_Y(y):=\sup_{\theta\in\mathbb{R}}\left\{\theta y-\frac{\theta^2}{2}\right\}\ (\mbox{for all}\ y\in\mathbb{R}),$$
which coincides with the rate function $J_Y(y)=\frac{y^2}{2}$.

Now we recall that $\sqrt{a_n}\frac{Y_n-\mu}{\sigma/\sqrt{n}}=\frac{nY_n-n\mu}{\sigma\sqrt{n/a_n}}$ and
$a_nnL(n)(Z_n-M)=z_n$ if and only if $Z_n=M+\frac{z_n}{a_nnL(n)}$; then, for $n$ large enough, we have
$$P\left(\sqrt{a_n}\frac{Y_n-\mu}{\sigma/\sqrt{n}}\in\cdot\Big|a_nnL(n)(Z_n-M)=z_n\right)
=P\left(\frac{M+\frac{z_n}{a_nnL(n)}+S_{n-1}^{(z_n)}-n\mu}{\sigma\sqrt{n/a_n}}\in\cdot\right),$$
where $S_{n-1}^{(z_n)}$ is the sum of $n-1$ i.i.d. random variables $W_1^{(z_n)},\ldots,W_{n-1}^{(z_n)}$ such that
$$\log\mathbb{E}\left[e^{\theta W_1^{(z_n)}}\right]=\kappa_{Y|Z}\left(\theta\Big|M+\frac{z_n}{a_nnL(n)}\right)
\ (\mbox{for all}\ \theta\in\mathbb{R}).$$
Thus we get
\begin{multline*}
\log\mathbb{E}\left[\exp\left(\frac{\theta}{a_n}\sqrt{a_n}\frac{Y_n-\mu}{\sigma/\sqrt{n}}\right)\Big|a_nnL(n)(Z_n-M)=z_n\right]\\
=(n-1)\kappa_{Y|Z}\left(\frac{\theta}{\sigma\sqrt{a_nn}}\Big|M+\frac{z_n}{a_nnL(n)}\right)
+\theta\,\frac{M+\frac{z_n}{a_nnL(n)}-n\mu}{\sigma\sqrt{a_nn}},
\end{multline*}
where, for a suitable remainder $o\left(\frac{1}{an_n}\right)$ (as $n\to\infty$),
\begin{multline*}
\kappa_{Y|Z}\left(\frac{\theta}{{\sigma\sqrt{a_nn}}}\Big|M+\frac{z_n}{a_nnL(n)}\right)=
\partial_\theta \kappa_{Y|Z}(0|M)\frac{\theta}{\sigma\sqrt{a_nn}}
+\partial_z \kappa_{Y|Z}(0|M)\frac{z_n}{a_nnL(n)}\\
+\frac{1}{2}\partial_{\theta\theta}^2\kappa_{Y|Z}(0|M)\frac{\theta^2}{\sigma^2a_nn}
+\frac{1}{2}\partial_{zz}^2\kappa_{Y|Z}(0|M)\frac{z_n^2}{a_n^2n^2L^2(n)}\\
+\partial_{\theta z}^2\kappa_{Y|Z}(0|M)\frac{\theta}{\sigma\sqrt{a_nn}}\frac{z_n}{a_nnL(n)}+o\left(\frac{1}{a_nn}\right);
\end{multline*}
moreover we have
$$\partial_\theta \kappa_{Y|Z}(0|M)=\mu,\ \partial_{\theta\theta}^2\kappa_{Y|Z}(0|M)=\sigma^2$$
and (we recall that $F(M)=1$, $\int_m^Mf(w)dw=1$ and $f(M)=F^\prime(M-)$ is finite and positive)
$$\partial_z \kappa_{Y|Z}(0|M)=\left.\frac{F(z)}{\int_m^ze^{\theta w}f(w)dw}
\frac{e^{\theta z}f(z)F(z)-f(z)\int_m^ze^{\theta w}f(w)dw}{F^2(z)}\right|_{(\theta,z)=(0,M)}=0.$$
Then we get the limit in \eqref{eq:GETlimit-conditionals-NCMD} noting that
\begin{multline*}
\frac{1}{1/a_n}\log\mathbb{E}\left[\exp\left(\frac{\theta}{a_n}\sqrt{a_n}\frac{Y_n-\mu}{\sigma/\sqrt{n}}\right)\Big|a_nnL(n)(Z_n-M)=z_n\right]\\
=a_n(n-1)\left\{\mu\frac{\theta}{\sigma\sqrt{a_nn}}+\frac{\sigma^2}{2}\frac{\theta^2}{\sigma^2a_nn}
+\frac{1}{2}\partial_{zz}^2\kappa_{Y|Z}(0|M)\frac{z_n^2}{a_n^2n^2L^2(n)}\right.\\
\left.+\partial_{\theta z}^2\kappa_{Y|Z}(0|M)\frac{\theta}{\sigma\sqrt{a_nn}}\frac{z_n}{a_nnL(n)}+o\left(\frac{1}{a_nn}\right)\right\}
+a_n\theta\,\frac{M+\frac{z_n}{a_nnL(n)}-n\mu}{\sigma\sqrt{a_nn}}\\
=\frac{\theta}{\sigma\sqrt{a_nn}}\left(a_n(n-1)\mu+\partial_{\theta z}^2\kappa_{Y|Z}(0|M)\frac{z_n(n-1)}{nL(n)}
+a_nM+\frac{z_n}{nL(n)}-a_nn\mu\right)\\
+\frac{\theta^2(n-1)}{2n}+\frac{a_n(n-1)}{2}\partial_{zz}^2\kappa_{Y|Z}(0|M)\frac{z_n^2}{a_n^2n^2L^2(n)}+a_n(n-1)o\left(\frac{1}{a_nn}\right)
\to\frac{\theta^2}{2}
\end{multline*}
(for each fixed $\theta\in\mathbb{R}$).

\paragraph{The proof of the goodness of the rate function $J_{Y,Z}$.}
Here we have to check that, for every $\eta\geq 0$, every closed level set of $J_{Y,Z}$ is compact. This can be done noting that, 
for every $\eta\geq 0$, we have
\begin{multline*}
\{(y,z)\in\mathbb{R}\times(-\infty,0]:J_{Y,Z}(y,z)\leq\eta\}=\{(y,z)\in\mathbb{R}\times(-\infty,0]:J_Y(y)+J_Z(z)\leq\eta\}\\
\subset\{y\in\mathbb{R}:J_Y(y)\leq\eta\}\times\{z\in(-\infty,0]:J_Z(z)\leq\eta\},
\end{multline*}
where both $\{y\in\mathbb{R}:J_Y(y)\leq\eta\}$ and $\{z\in(-\infty,0]:J_Z(z)\leq\eta\}$ are compact sets; so every level set is 
compact because it is a subset of a compact set.
\end{proof}

\begin{remark}\label{rem:asympt-independence}
The rate function $J_{Y,Z}(y,z)$ in Proposition \ref{prop:sum-maxima-NCMD} can be expressed as a sum of two functions which 
depend on $y$ and $z$ only, i.e. the \emph{marginal rate functions} $J_Y(y)$ and $J_Z(z)$ that appear in the proof of that 
proposition. This is not surprising by the asymptotic independence stated in Proposition \ref{prop:ChowTeugels-consequence}.
\end{remark}

\section{A modification of Proposition \ref{prop:LDP-sum-maxima} when $M$ is not finite}\label{sec:added}
In this section we prove Proposition \ref{prop:LDP-sum-maxima-added}, i.e.  a suitable modification of Proposition 
\ref{prop:LDP-sum-maxima} with $\{P((Y_n,Z_n/h_n)\in\cdot):n\geq 1\}$ in place of $\{P((Y_n,Z_n)\in\cdot):n\geq 1\}$, 
for some $h_n$ such that $h_n\to\infty$; actually we consider some different hypotheses and, in particular, $M=\infty$.
In order to do that we refer to Proposition 3.1 in \cite{GiulianoMacciCSTM} (in place of Proposition 4.1 in \cite{GiulianoMacciCSTM};
we mean the part of the proof of Proposition \ref{prop:LDP-sum-maxima} in which we check that \textbf{(C1)} holds). We start with 
the following useful lemma.

\begin{lemma}\label{lem:two-speed-functions}
	Let $\{\pi_n\}_n$ be a sequence of probability measures (on some Polish space) that satisfies the LDP with speed $s_n$ and good rate 
	function $I$, which uniquely vanishes at some $r_0$. Moreover let $t_n$ be another speed function such that $\frac{s_n}{t_n}\to\infty$. 
	Then $\{\pi_n\}_n$ satisfies the LDP with speed $t_n$ and good rate function $\Delta(\cdot;r_0)$ defined by
	\begin{equation}
		\Delta(\cdot;r_0):=\left\{\begin{array}{ll}
			0&\ \mbox{if}\ r=r_0\\
			\infty&\ \mbox{if}\ r\neq r_0.
		\end{array}\right.
	\end{equation}
\end{lemma}

\begin{proof}
	Firsty we can say that $\{\pi_n\}_n$ is exponentially tight with respect to $s_n$ (this follows from the LDP of the 
	sequence $\{\pi_n\}_n$ with speed $s_n$ and good rate function $I$, and Lemma 2.6 in \cite{LynchSethuraman}). Then 
	$\{\pi_n\}_n$ is also exponentially tight with respect to $t_n$; indeed, if for every $b>0$ there exists a compact set $K_b$
	such that
	$$\pi_n(K_b^c)\leq ae^{-s_n b}\quad \mbox{eventually}$$
	for some $a>0$, then we have the same estimate with $t_n$ in place $s_n$ because $e^{-s_n}\leq e^{-t_n}$. So there exists at least 
	a subsequence of $\{\pi_n\}_n$ which satisfies the LDP with speed $t_n$ (see e.g. Theorem (P) in \cite{Puhalskii}). We complete the 
	proof showing that, for every subsequence of $\{\pi_n\}_n$ (which we still call $\{\pi_n\}_n$) that satisfies the LDP with speed 
	$t_n$, the governing rate function is $\Delta(\cdot;r_0)$. Here, as we did in Section \ref{sec:preliminaries}, we consider the 
	notation $B_R(r)$ for the open ball centered at $r$ and with radius $R$. The, by the hypotheses, we have
	$$-I(r)\leq\lim_{R\to 0}\liminf_{n\to\infty}\frac{1}{s_n}\log\pi_n(B_R(r))
	\leq\lim_{R\to 0}\limsup_{n\to\infty}\frac{1}{s_n}\log\pi_n(B_R(r))\leq-I(r)$$
	for every $r$ in the Polish space; our aim is to get the same estimate (up to a subsequence) with $t_n$ in place of $s_n$ 
	and $\Delta(\cdot;r_0)$ in place of $I$.
	
	We start with the case $r=r_0$. Then we trivially have
	$$\limsup_{n\to\infty}\frac{1}{t_n}\log\pi_n(B_R(r))\leq 0=-\Delta(r_0;r_0),$$
	whence we obtain
	$$\lim_{R\to 0}\limsup_{n\to\infty}\frac{1}{t_n}\log\pi_n(B_R(r))\leq-\Delta(r_0;r_0).$$
	Moreover, for every $R>0$, we have $\pi_n(B_R(r))\to 1$; this yields
	$$\lim_{n\to\infty}\frac{1}{t_n}\log\pi_n(B_R(r))=0=-\Delta(r_0;r_0),$$
	whence we obtain
	$$\lim_{R\to 0}\liminf_{n\to\infty}\frac{1}{t_n}\log\pi_n(B_R(r))=-\Delta(r_0;r_0).$$
	Thus the desired bounds for $r=r_0$ are proved. Now the case $r\neq r_0$. Then we trivially have
	$$\liminf_{n\to\infty}\frac{1}{t_n}\log\pi_n(B_R(r))\geq-\infty=-\Delta(r;r_0),$$
	whence we obtain
	$$\lim_{R\to 0}\liminf_{n\to\infty}\frac{1}{t_n}\log\pi_n(B_R(r))\geq-\Delta(r;r_0).$$
	Moreover we can find $\rho>0$ small enough to have $I(\overline{B_\rho(r)}):=\inf\{I(y):y\in\overline{B_\rho(r)}\}>0$
	(thus $r_0\notin\overline{B_\rho(r)}$). Then
	$$\limsup_{n\to\infty}\frac{1}{t_n}\log\pi_n(B_R(r))
	\leq\limsup_{n\to\infty}\frac{s_n}{t_n}\frac{1}{s_n}\log\pi_n(\overline{B_\rho(r)})
	\leq-\infty=-\Delta(r;r_0)$$
	(because $\frac{s_n}{t_n}\to\infty$ and $\limsup_{n\to\infty}\frac{1}{s_n}\log\pi_n(\overline{B_\rho(r)})\leq-I(\overline{B_\rho(r)})$);
	so, by the monotonicity with respect to $\rho$, we get
	$$\lim_{R\to 0}\limsup_{n\to\infty}\frac{1}{t_n}\log\pi_n(B_R(r))
	\leq\limsup_{n\to\infty}\frac{1}{t_n}\log\pi_n(B_\rho(r))\leq-\Delta(r;r_0).$$
	Thus the desired bounds for $r\neq r_0$ are proved, and this completes the proof.
\end{proof}

Now we are able to prove Proposition \ref{prop:LDP-sum-maxima-added}. In particular we consider the notation in Assumption 
\ref{assumpt:LD}, and, again, we use the notation $\mu$ for the mean of the i.i.d. random variables $\{W_n:n\geq 1\}$.

\begin{proposition}\label{prop:LDP-sum-maxima-added}
	Let $\{W_n:n\geq 1\}$ be i.i.d. random variables with common continuous distribution function $F$ such that 
	$\kappa_Y(\theta):=\log\mathbb{E}[e^{\theta W_1}]$ is finite in a neighbourhood of $\theta=0$. Assume that
	$M=\infty$. We set $\mathcal{H}(x)=-\log(1-F(x))$. Moreover, let $h_n$ be such that $1-F(h_n)=\frac{1}{n}$, or equivalently
	$\mathcal{H}(h_n)=\log n$. We also assume that $\mathcal{H}$ is a regularly varying function at $\infty$ of index $\alpha>0$, i.e.
	$$\lim_{y\to\infty}\frac{\mathcal{H}(xy)}{\mathcal{H}(y)}=x^\alpha\quad \mbox{for all}\ x>0.$$
	Then $\{P((Y_n,Z_n/h_n)\in\cdot):n\geq 1\}$ satisfies the LDP with speed $\log n$ and rate function $H_{Y,Z}$ defined by
    $$H_{Y,Z}(y,z):=\left\{\begin{array}{ll}
             H_Z(z)&\ \mbox{if}\ z\geq 1\ \mbox{and}\ y=\mu\\
    	     \infty&\ \mbox{otherwise},
    \end{array}\right.$$
    where $H_Z(z):=z^\alpha-1$.
\end{proposition}
\begin{proof}
	It is well-known that it is enough to prove the two following conditions:
	\begin{enumerate}
	    \item for all $(y,z)\in\mathbb{R}^2$
	    \begin{multline*}
	    	-H_{Y,Z}(y,z)\leq\lim_{R\to 0}\liminf_{n\to\infty}\frac{1}{\log n}\log P((Y_n,Z_n/h_n)\in(y-R,y+R)\times(z-R,z+R))\\
	    	\leq\lim_{R\to 0}\limsup_{n\to\infty}\frac{1}{\log n}\log P((Y_n,Z_n/h_n)\in(y-R,y+R)\times(z-R,z+R))\leq-H_{Y,Z}(y,z);
	    \end{multline*}
	    \item $\{P((Y_n,Z_n/h_n)\in\cdot):n\geq 1\}$ is exponentially tight with respect to the speed $\log n$.
	\end{enumerate}
	 For the first condition we start with two trivial cases $z<1$ and $y\neq\mu$, and it is enough to
	 check the upper bound. If $z<1$ we have
	 \begin{equation}\label{eq:trivial-bound}
	 	P((Y_n,Z_n/h_n)\in(y-R,y+R)\times(z-R,z+R))\leq P(Z_n/h_n\in(z-R,z+R))
	 \end{equation}
	 and, for $R>0$ small enough,
	 $$\limsup_{n\to\infty}\frac{1}{\log n}\log P(Z_n/h_n\in(z-R,z+R))=-\infty$$
	 by the LDP in Proposition 3.1 in \cite{GiulianoMacciCSTM}. If $y\neq\mu$ we have
	 $$P((Y_n,Z_n/h_n)\in(y-R,y+R)\times(z-R,z+R))
	 \leq P(Y_n\in(y-R,y+R))$$
	 and, for $R>0$ small enough,
	 $$\limsup_{n\to\infty}\frac{1}{\log n}\log P(Y_n\in(y-R,y+R))=-\infty$$
	 by the LDP of $\{Y_n:n\geq 1\}$ with speed $\log n$ with rate function $\Delta(\cdot;\mu)$ in Lemma \ref{lem:two-speed-functions};
	 this LDP is a consequence of Lemma \ref{lem:two-speed-functions} together with Cramér Theorem (see e.g. Theorem 2.2.3 in \cite{DemboZeitouni})
	 with $I=\kappa_Y^*$ (where $\kappa_Y^*$ defined by
	 $$\kappa_Y^*(y):=\sup_{\theta\in\mathbb{R}}\{\theta y-\kappa_Y(\theta)\}$$
	 which uniquely vanishes at $y=\mu$), $s_n=n$ and $t_n=\log n$.
	 
	 So we conclude the proof of the first condition by taking $z\geq 1$ and $y=\mu$. The upper bound can be proved as we did before 
	 for the case $z<1$; indeed, by \eqref{eq:trivial-bound} and by the LDP in Proposition 3.1 in \cite{GiulianoMacciCSTM}, we have
	 \begin{multline*}
	 	\lim_{R\to 0}\limsup_{n\to\infty}\frac{1}{\log n}\log P((Y_n,Z_n/h_n)\in(y-R,y+R)\times(z-R,z+R))\\
	 	\leq\lim_{R\to 0}\limsup_{n\to\infty}\frac{1}{\log n}\log P(Z_n/h_n\in(z-R,z+R))\leq-H_Z(z).
	 \end{multline*}
     For the lower bound we take into account that
     \begin{multline*}
     	P((Y_n,Z_n/h_n)\in(y-R,y+R)\times(z-R,z+R))\\
     	=P(Z_n/h_n\in(z-R,z+R))-P((Y_n,Z_n/h_n)\in(y-R,y+R)^c\times(z-R,z+R)),
     \end{multline*}
     and we get
     $$\lim_{R\to 0}\liminf_{n\to\infty}\frac{1}{\log n}\log P((Y_n,Z_n/h_n)\in(y-R,y+R)\times(z-R,z+R))\geq-H_Z(z)$$
     by applying Lemma 19 in \cite{GaneshTorrisi}. In order to do that we remark that
     $$\liminf_{n\to\infty}\frac{1}{\log n}\log P(Z_n/h_n\in(z-R,z+R))\geq-H_Z(z)$$
     by the LDP in Proposition 3.1 in \cite{GiulianoMacciCSTM}, and
     \begin{multline*}
         \limsup_{n\to\infty}\frac{1}{\log n}\log P((Y_n,Z_n/h_n)\in(y-R,y+R)^c\times(z-R,z+R))\\
         \leq\limsup_{n\to\infty}\frac{1}{\log n}\log P(Y_n\in(y-R,y+R)^c)\leq-\inf_{s\in(y-R,y+R)^c}\Delta(s,\mu)=-\infty
     \end{multline*}
     (here we take into account the LDP of $\{Y_n:n\geq 1\}$ with speed $\log n$ stated above). Then Lemma 19 in \cite{GaneshTorrisi}
     yields
     $$\liminf_{n\to\infty}\frac{1}{\log n}\log P((Y_n,Z_n/h_n)\in(y-R,y+R)\times(z-R,z+R))\geq -H_Z(z),$$
     and we easily get the desired lower bound.
     
     We conclude with the second condition, i.e. the exponential tightness. By Lemma 2.6 in \cite{LynchSethuraman} the marginal sequences
     are exponentially tight; thus, for all $b>0$, there exist two compact sets $K_b^{(1)}$ and $K_b^{(2)}$ such that
     $$P(Y_n\notin K_b^{(1)})\leq a_1e^{-b\log n}\ \mbox{and}\ P(Z_n/h_n\notin K_b^{(2)})\leq a_2e^{-b\log n}\quad\mbox{eventually},$$
     for some $a_1,a_2>0$. Then, since $K_z^{(1)}\times K_z^{(2)}$ is a compact set, we conclude the proof noting that
     $$P((Y_n,Z_n/h_n)\notin K_b^{(1)}\times K_b^{(2)})\leq P(Y_n\notin K_b^{(1)})+P(Z_n/h_n\notin K_b^{(2)})\leq (a_1+a_2)e^{-b\log n}\quad\mbox{eventually}.$$
\end{proof}

We conclude noting that, as it happens for the rate function $J_{Y,Z}(y,z)$ in Proposition \ref{prop:sum-maxima-NCMD} (see Remark 
\ref{rem:asympt-independence}), we have an asymptotic independence interpretation for the rate function $H_{Y,Z}(y,z)$ in Proposition 
\ref{prop:LDP-sum-maxima-added}.

\begin{remark}\label{rem:asympt-independence-added}
	The rate function $H_{Y,Z}(y,z)$ in Proposition \ref{prop:LDP-sum-maxima-added} can be expressed as a sum of two functions which 
	depend on $y$ and $z$ only, i.e. the \emph{marginal rate functions} $\Delta(y;\mu)$ and $H_Z(z)$ that appear in the proof of that 
	proposition.
\end{remark}

\section{MD for sums of minima of i.i.d. exponential random variables}\label{sec:result2}
We start with the following assumption.

\begin{assumption}\label{assumpt:minima-iid-exp}
Let $\{W_n:n\geq 1\}$ be a sequence of i.i.d. real random variables with exponential distribution; more precisely their common
distribution function $F$ is defined by
$$F(x):=1-e^{-\lambda x}\ \mbox{for all}\ x\geq 0.$$
Moreover let $\{X_n:n\geq 1\}$ be the sequence of random variables defined by
$$X_n:=\frac{\sum_{k=1}^n\min\{W_1,\ldots,W_k\}}{\log n}\ \mbox{for all}\ n\geq 2.$$
\end{assumption}

Now we recall two results. The first one provides the reference LDP, namely the LDP which governs the convergence of $X_n$ to 
$\frac{1}{\lambda}$ (as $n\to\infty$); indeed the rate function $I_X$ in the next proposition uniquely vanishes at 
$x=\frac{1}{\lambda}$.

\begin{proposition}\label{prop:GiulianoMacciESAIM}
Assume that Assumption \ref{assumpt:minima-iid-exp} holds. Then $\{P(X_n\in\cdot):n\geq 2\}$
satisfies the LDP with speed $\log n$ and rate function $I_X$ defined by
$$I_X(x):=\left\{\begin{array}{ll}
(\sqrt{\lambda x}-1)^2&\ \mbox{if}\ x\geq 0\\
\infty&\ \mbox{if}\ x<0.
\end{array}\right.$$
\end{proposition}
\begin{proof}
See Proposition 5.2 in \cite{GiulianoMacciESAIM}.
\end{proof}

The second result concerns the following weak convergence to a centered Gaussian distribution.

\begin{proposition}\label{prop:Hoglund-consequence}
Assume that Assumption \ref{assumpt:minima-iid-exp} holds. Then $(X_n-\frac{1}{\lambda})\sqrt{\log n}$
converges weakly (as $n\to\infty$) to the centered Gaussian distribution with variance 
$\sigma^2=\frac{2}{\lambda^2}$.
\end{proposition}
\begin{proof}
The random variables $(X_n-\frac{1}{\lambda})\sqrt{\frac{\lambda^2}{2}\log n}$ converge weakly to the standard Gaussian 
di\-stribution by Theorem in \cite{Hoglund}; indeed the distribution function $F$ in Assumption \ref{assumpt:minima-iid-exp}
satisfies the condition $\int_0^1|F(x)-x/b|x^{-2}dx<\infty$ (required in \cite{Hoglund}) if and only if $b=\frac{1}{\lambda}$.
Then we can immediately conclude with the desired weak convergence.
\end{proof}

The aim of this section is to prove Proposition \ref{prop:MD} which provides a class of LDPs that fills the 
gap between the convergence of $\{X_n:n\geq 1\}$ to $\frac{1}{\lambda}$ (governed by the LDP in Proposition 
\ref{prop:GiulianoMacciESAIM} with speed $v_n=\log n$), and the weak convergence in Proposition 
\ref{prop:Hoglund-consequence}. Then we get a (\emph{central}) moderate deviation result because the weak convergence in 
Proposition \ref{prop:Hoglund-consequence} is towards a Gaussian distribution. We also remark that, as it typically happens,
we have $I_X^{\prime\prime}(\frac{1}{\lambda})=\frac{1}{\sigma^2}$ (where $\sigma^2=\frac{2}{\lambda^2}$ as in
Proposition \ref{prop:Hoglund-consequence}); this equality can be checked with some easy computations, and we 
omit the details.

\begin{proposition}\label{prop:MD}
Assume that Assumption \ref{assumpt:minima-iid-exp} holds. Then, for every sequence of positive numbers 
$\{a_n:n\geq 1\}$ such that \eqref{eq:MD-conditions} holds with $v_n=\log n$, the sequence 
$\left\{P\left(\left(X_n-\frac{1}{\lambda}\right)\sqrt{a_n\log n}\in\cdot\right):n\geq 2\right\}$ 
satisfies the LDP with speed $1/a_n$ and rate function $J_X$ defined by $J_X(x)=\frac{x^2}{2\sigma^2}$, where 
$\sigma^2=\frac{2}{\lambda^2}$ as in Proposition \ref{prop:Hoglund-consequence}.
\end{proposition}
\begin{proof}
We apply G\"artner Ellis Theorem, i.e. Proposition \ref{prop:GET}. Indeed we show that
\begin{equation}\label{eq:GETlimit-MD}
\lim_{n\to\infty}\frac{1}{1/a_n}\log\mathbb{E}\left[\exp\left(\frac{\theta}{a_n}\left(X_n-\frac{1}{\lambda}\right)
\sqrt{a_n\log n}\right)\right]=\underbrace{\frac{\sigma^2\theta^2}{2}}_{=\theta^2/\lambda^2}\ (\mbox{for all}\ \theta\in\mathbb{R})
\end{equation}
and therefore we get the desired LDP with rate function $J_X$ defined by
$$J_X(x):=\sup_{\theta\in\mathbb{R}}\left\{\theta x-\frac{\theta^2}{\lambda^2}\right\}\ (\mbox{for all}\ x\in\mathbb{R}),$$
which coincides with the rate function $J_X$ in the statement.
	
We use a known expression for the moment generating function of $\sum_{k=1}^n\min\{W_1,\ldots,W_k\}$ (see e.g. eq. (3.5) in 
\cite{GhoshBabuMukhopadhyay}):
\begin{multline*}
\frac{1}{1/a_n}\log\mathbb{E}\left[\exp\left(\frac{\theta}{a_n}\left(X_n-\frac{1}{\lambda}\right)\sqrt{a_n\log n}\right)\right]
=a_n\left(-\frac{\theta\sqrt{a_n\log n}}{\lambda a_n}+\log\mathbb{E}\left[\exp\left(\frac{\theta\sqrt{a_n\log n}}{a_n}X_n\right)\right]\right)\\
=-\frac{\theta\sqrt{a_n\log n}}{\lambda}+a_n\log\mathbb{E}\left[\exp\left(\frac{\theta\sum_{k=1}^n\min\{W_1,\ldots,W_k\}}{\sqrt{a_n\log n}}\right)\right]\\
=\left\{\begin{array}{ll}
-\frac{\theta\sqrt{a_n\log n}}{\lambda}+a_n\sum_{k=1}^n\log\left(1+\frac{\frac{\theta}{\lambda\sqrt{a_n\log n}}}
{k\left(1-\frac{\theta}{\lambda\sqrt{a_n\log n}}\right)}\right)&\ \mbox{if}\ \frac{\theta}{\lambda\sqrt{a_n\log n}}<1\\
\infty&\ \mbox{otherwise};
\end{array}\right.
\end{multline*}
then, for each fixed $\theta\in\mathbb{R}$, we can take $n$ large enough in order to have $\frac{\theta}{\lambda\sqrt{a_n\log n}}<1$
(since $a_n\log n\to\infty$, as $n\to\infty$).
	
Moreover we remark that
\begin{equation}\label{eq:local-estimate-degree2}
\mbox{for all}\ v>\frac{1}{2},\ \mbox{there exists}\ \delta>0\
\mbox{such that}\ \log(1+x)\geq x-vx^2\ \mbox{for all}\ |x|<\delta
\end{equation}
(this can be proved by checking that the function $g$ defined by $g(x):=\log(1+x)-(x-vx^2)$ has a local minimum at $x=0$);
so, for $\delta>0$ as in \eqref{eq:local-estimate-degree2}, we take $n$ large enough in order to have 
$\left|\frac{\theta}{\lambda\sqrt{a_n\log n}}/\left(1-\frac{\theta}{\lambda\sqrt{a_n\log n}}\right)\right|<\delta$.
	
Finally we set
\begin{multline*}
b_n:=-\frac{\theta\sqrt{a_n\log n}}{\lambda}+\frac{\frac{a_n\theta}{\lambda\sqrt{a_n\log n}}}
{1-\frac{\theta}{\lambda\sqrt{a_n\log n}}}\sum_{k=1}^n\frac{1}{k}\\
=\frac{-\frac{\theta\sqrt{a_n\log n}}{\lambda}+\frac{\theta^2}{\lambda^2}+\frac{a_n\theta}{\lambda\sqrt{a_n\log n}}\sum_{k=1}^n\frac{1}{k}}
{1-\frac{\theta}{\lambda\sqrt{a_n\log n}}}
=\frac{\frac{\theta\sqrt{a_n\log n}}{\lambda}\left(-1+\frac{\sum_{k=1}^n1/k}{\log n}\right)+\frac{\theta^2}{\lambda^2}}
{1-\frac{\theta}{\lambda\sqrt{a_n\log n}}},
\end{multline*}
and, for $n$ large enough, we have
$$b_n-v\frac{\frac{\theta^2}{\lambda^2\log n}}
{\left(1-\frac{\theta}{\lambda\sqrt{a_n\log n}}\right)^2}\sum_{k=1}^n\frac{1}{k^2}\leq-\frac{\theta\sqrt{a_n\log n}}{\lambda}+a_n\sum_{k=1}^n\log\left(1+\frac{\frac{\theta}{\lambda\sqrt{a_n\log n}}}
{k\left(1-\frac{\theta}{\lambda\sqrt{a_n\log n}}\right)}\right)\leq b_n,$$
by using \eqref{eq:local-estimate-degree2} with 
$x=\frac{\theta}{\lambda\sqrt{a_n\log n}}/\left(1-\frac{\theta}{\lambda\sqrt{a_n\log n}}\right)$, and by the 
well-known inequality $\log (1+y)\leq y$, for every $y>-1$.
	
So the desired condition \eqref{eq:GETlimit-MD} holds since
\begin{equation}\label{eq:final-limits}
\lim_{n\to\infty}b_n=\frac{\theta^2}{\lambda^2}\quad \mbox{and}\quad \lim_{n\to\infty}\frac{\frac{\theta^2}{\lambda^2\log n}}
{\left(1-\frac{\theta}{\lambda\sqrt{a_n\log n}}\right)^2}\sum_{k=1}^n\frac{1}{k^2}=0.
\end{equation}
Indeed the first limit in \eqref{eq:final-limits} holds by \eqref{eq:MD-conditions} with $v_n=\log n$ (which yields 
$a_n\to 0$ and $a_n\log n\to\infty$), and by
$$\lim_{n\to\infty}\sqrt{\log n}\left(-1+\frac{\sum_{k=1}^n1/k}{\log n}\right)=0;$$
the second limit in \eqref{eq:final-limits} trivially holds by taking into account $a_n\log n\to\infty$ and 
$\sum_{k=1}^\infty\frac{1}{k^2}<\infty$.
\end{proof}

In order to make the paper more self-contained we remark that the limit in \eqref{eq:GETlimit-MD} with 
$a_n=1$ yields the weak convergence in Proposition \ref{prop:Hoglund-consequence}.

\paragraph{Funding.}
This work has been partially supported by MIUR Excellence Department Project awarded to the Department of Mathematics, 
University of Rome Tor Vergata (CUP E83C18000100006 and CUP E83C23000330006), by University of Rome Tor Vergata (project 
"Asymptotic Methods in Probability" (CUP E89C20000680005) and project "Asymptotic Properties in Probability" 
(CUP E83C22001780005)) and by Indam-GNAMPA.

\paragraph{Acknowledgements.}
The authors thank two referees for some useful comments and Professor Clive W. Anderson for some discussion on the 
content of reference \cite{ChowTeugels}.


\begin{thebibliography}{spc}
\bibitem{AndersonTurkman}
C.W. Anderson, K.F. Turkman (1991) The joint limiting distribution of sums and maxima of stationary sequences. 
J. Appl. Probab. 28, no. 1, 33--44.
\bibitem{ArendarczykKozubowskiPanorska}
M. Arendarczyk, T.J. Kozubowski, A.K. Panorska (2018) The joint distribution of the sum and maximum of 
dependent Pareto risks. J. Multivariate Anal. 167, 136--156.
\bibitem{Chaganty}
N.R. Chaganty (1997) Large deviations for joint distributions and statistical applications. Sankhy\={a}  A
59, no. 2, 147--166.
\bibitem{ChowTeugels}
T.L. Chow, J.L. Teugels (1979) The sum and the maximum of i.i.d. random variables.
Proceedings of the Second Prague Symposium on Asymptotic Statistics (Hradec Kr\'alov\'e,
1978), pp. 81--92, North-Holland, Amsterdam-New York.
\bibitem{DemboZeitouni}
A. Dembo, O. Zeitouni (1998) Large Deviations Techniques and Applications, 2nd  edn. Springer.
\bibitem{EmbrechtsKluppelbergMikosch}
P. Embrechts, C. Klüppelberg, T. Mikosch (1997) Modelling Extremal Events. Springer-Verlag.
\bibitem{GaneshTorrisi}
A. Ganesh, G.L. Torrisi (2008) Large deviations of the interference in a wireless communication model.
IEEE Trans. Inform. Theory 54, 3505--3517.
\bibitem{GhoshBabuMukhopadhyay}
M. Ghosh, G.J. Babu, N. Mukhopadhyay (1975) Almost sure convergence of sums of maxima and minima of positive 
random variables. Z. Wahrsch. Verw. Gebiete 33, 49--54.
\bibitem{GiulianoMacciCSTM}
R. Giuliano, C. Macci (2014) Large deviation principles for sequences of maxima and minima. Comm. Statist. 
Theory Methods 43, no. 6, 1077--1098.
\bibitem{GiulianoMacciESAIM}
R. Giuliano, C. Macci (2015) Asymptotic results for weighted means of random variables which converge to 
a Dickman distribution, and some number theoretical applications. ESAIM Probab. Stat. 19, 395--413.
\bibitem{GiulianoMacciMSTA}
R. Giuliano, C. Macci (2023) Some examples of noncentral moderate deviations for sequences 
of real random variables. Mod. Stoch. Theory Appl. 10, no. 2, (2023) 111--144.
\bibitem{Hoglund}
T. H\"oglund (1972) Asymptotic normality of sums of minima of random variables. Ann. Math. Statist. 43,
351--353.
\bibitem{Hsing}
T. Hsing (1995) A note on the asymptotic independence of the sum and maximum of strongly mixing stationary
random variables. Ann. Probab. 23, 938--947.
\bibitem{Kratz}
M. Kratz (2014) Normex, a new method for evaluating the distribution of aggregated heavy tailed risks.
Extremes 17, 661--691.
\bibitem{KratzProkopenko}
M. Kratz, E. Prokopenko (2023+) Multi-normex distributions for the sum of random vectors. Rates of convergence.
Extremes, to appear.
\bibitem{Krizmanic}
D. Krizmanić (2020) On joint weak convergence of partial sum and maxima processes. Stochastics 92, no. 6, 
876--899.
\bibitem{LeonenkoMacciPacchiarotti}
N. Leonenko, C. Macci, B. Pacchiarotti (2021) Large deviations for a class of tempered subordinators and their 
inverse processes. Proc. Roy. Soc. Edinburgh Sect. A 151, no. 6, 2030--2050.
\bibitem{LynchSethuraman}
J. Lynch, J. Sethuraman (1987) Large deviations for processes with independent increments. Ann. Probab. 15, 
610--627.
\bibitem{Muller}
U.K. M\"uller (2019) Refining the central limit theorem approximation via extreme value theory. Statist. 
Probab. Lett. 155, Article ID 108564, 7 pages.
\bibitem{Puhalskii}
A. Puhalskii (1991) On functional principle of large deviations. New Trends in Probability and Statistics, 
vol. 1, pp. 198--218.
\bibitem{QeadanKozubowskiPanorska}
F. Qeadan, T.J. Kozubowski, A.K. Panorska (2012) The joint distribution of the sum and the maximum of IID 
exponential random variables. Comm. Statist. Theory Methods 41, no. 3, 544--569.
\end{thebibliography}
\end{document}